\documentclass[a4paper]{article}

\usepackage{amssymb,amsmath, amsfonts}
\usepackage{authblk}
\usepackage{pst-grad} 

\def\str{\rightarrow}
\def\mj{\mbox{\bf 1}}
\def\m1{\begin{picture}(14,10)(0,2)
\put(7,0){\line(-1,2){5}} \put(7,0){\line(1,2){5}}
\put(7,0){\circle*{1.5}} \put(2,10){\circle*{1.5}}
\put(12,10){\circle*{1.5}}\end{picture}}
\def\e1{\begin{picture}(10,10)(0,2) \put(5,0){\circle*{1.5}}\end{picture}}
\def\p1{\begin{picture}(10,10)(0,2) \put(5,10){\circle*{1.5}}\end{picture}}
\def\pl{\!+\!}
\def\mn{\!-\!}
\def\dkz{\noindent{\sc Proof. }}
\def\qed{\hfill $\dashv$}
\def\top{\mbox{\it{Top}}}

\def\mm2{\begin{picture}(34,10)(-10,2)
\put(7,0){\line(-1,2){5}} \put(7,0){\line(1,2){5}}
\put(-6,0){\line(0,1){10}} \put(20,0){\line(0,1){10}}
\put(7,0){\circle*{1.5}} \put(2,10){\circle*{1.5}}
\put(12,10){\circle*{1.5}} \put(-6,0){\circle*{1.5}}
\put(20,0){\circle*{1.5}} \put(-6,10){\circle*{1.5}}
\put(20,10){\circle*{1.5}}\end{picture}}

\def\ee2{\begin{picture}(30,10)(-10,2) \put(5,0){\circle*{1.5}} \put(-3,0){\circle*{1.5}} \put(13,0){\circle*{1.5}}
\put(-3,10){\circle*{1.5}} \put(13,10){\circle*{1.5}}
\put(-3,0){\line(0,1){10}} \put(13,0){\line(0,1){10}}
\end{picture}}

\def\pp2{\begin{picture}(24,10)(0,2)
\put(7,0){\line(-1,2){5}} \put(7,0){\line(1,2){5}}
\put(20,0){\line(0,1){10}} \put(7,0){\circle*{1.5}}
\put(2,10){\circle*{1.5}} \put(12,10){\circle*{1.5}}
\put(20,0){\circle*{1.5}} \put(20,10){\circle*{1.5}}\end{picture}}

\def\ppp3{\begin{picture}(24,10)(-10,2)
\put(7,0){\line(-1,2){5}} \put(7,0){\line(1,2){5}}
\put(-6,0){\line(0,1){10}} \put(7,0){\circle*{1.5}}
\put(2,10){\circle*{1.5}} \put(12,10){\circle*{1.5}}
\put(-6,0){\circle*{1.5}} \put(-6,10){\circle*{1.5}}
\end{picture}}

\def\i3{\begin{picture}(34,10)(0,2)
\put(7,0){\line(-1,2){5}} \put(7,0){\line(1,2){5}}
\put(7,0){\line(0,1){10}} \put(20,0){\line(0,1){10}}
\put(30,0){\line(0,1){10}} \put(7,0){\circle*{1.5}}
\put(7,10){\circle*{1.5}} \put(2,10){\circle*{1.5}}
\put(12,10){\circle*{1.5}} \put(20,0){\circle*{1.5}}
\put(20,10){\circle*{1.5}} \put(30,0){\circle*{1.5}}
\put(30,10){\circle*{1.5}}\end{picture}}

\def\k1{\begin{picture}(34,10)(0,2)
\put(0,0){\line(0,1){10}} \put(8,0){\line(0,1){10}}
\put(16,0){\line(0,1){10}}\put(27,0){\line(-1,2){5}}
\put(27,0){\line(1,2){5}} \put(27,0){\circle*{1.5}}
\put(22,10){\circle*{1.5}} \put(32,10){\circle*{1.5}}
\put(0,0){\circle*{1.5}} \put(0,10){\circle*{1.5}}
\put(8,0){\circle*{1.5}} \put(8,10){\circle*{1.5}}
\put(16,0){\circle*{1.5}} \put(16,10){\circle*{1.5}}\end{picture}}

\def\kk2{\begin{picture}(36,10)(0,2)
\put(18,0){\line(0,1){10}} \put(26,0){\line(0,1){10}}
\put(34,0){\line(0,1){10}}\put(7,0){\line(-1,2){5}}
\put(7,0){\line(1,2){5}} \put(7,0){\circle*{1.5}}
\put(2,10){\circle*{1.5}} \put(12,10){\circle*{1.5}}
\put(18,0){\circle*{1.5}} \put(18,10){\circle*{1.5}}
\put(26,0){\circle*{1.5}} \put(26,10){\circle*{1.5}}
\put(34,0){\circle*{1.5}} \put(34,10){\circle*{1.5}}\end{picture}}

\allowdisplaybreaks

 \newtheorem{prop}{Proposition}[section]

\begin{document}

\title{Monoids, Segal's condition and bisimplicial spaces}
\author{Zoran Petri\'c}
\affil{Mathematical Institute, SANU,\\ Knez Mihailova 36, p.f.\
367,\\ 11001 Belgrade, Serbia\\ \texttt{zpetric@mi.sanu.ac.rs} }

\date{}
\maketitle

\vspace{-3ex}

\begin{abstract}
A characterization of simplicial objects in categories with finite
products obtained by the reduced bar construction is given. The
condition that characterizes such simplicial objects is a
strictification of Segal's condition guaranteeing that the loop
space of the geometric realization of a simplicial space $X$ and
the space $X_1$ are of the same homotopy type. A generalization of
Segal's result appropriate for bisimplicial spaces is given. This
generalization gives conditions guaranteing that the double loop
space of the geometric realization of a bisimplicial space $X$ and
the space $X_{11}$ are of the same homotopy type.
\end{abstract}

\vspace{.3cm}


\vspace{.5ex}

\noindent {\small {\it Keywords$\,$}: reduced bar construction,
simplicial space, loop space}

\vspace{.5ex}

\noindent {\small {\it Acknowledgements$\,$}: I am grateful to
Rade \v Zivaljevi\' c and Matija Ba\v si\' c for some useful
discussions during the CGTA Colloquium in Belgrade. This work was
supported by project ON174026 of the Ministry of Education,
Science, and Technological Development of the Republic of Serbia.
}

\section{Introduction}

This paper is based on the author talks delivered in 2014 at the
Fourth Mathematical Conference of the Republic of Srpska and at
the CGTA Colloquium of the Faculty of Mathematics in Belgrade. Its
first part gives a condition which is necessary and sufficient for
a simplicial object to be obtained by the reduced bar
construction. It turns out that this condition is a
strictification of Segal's condition guaranteeing that the loop
space of the geometric realization of a simplicial space $X$ and
the space $X_1$ are of the same homotopy type.

The second part of this paper is devoted to a generalization of
Segal's result. This generalization gives conditions guaranteing
that the double loop space of the geometric realization of a
bisimplicial space $X$ and the space $X_{11}$ are of the same
homotopy type. We refer to \cite{P14} for a complete
generalization of Segal's result.

\section{Monoids and the reduced bar construction}

A \emph{strict monoidal} category $(\cal M,\otimes,I)$ is a
category $\cal M$ with an associative bifunctor $ \otimes\!:{\cal
M}\times {\cal M}\str {\cal M}$,
\[(A\otimes B)\otimes C=A\otimes(B\otimes C)\quad {\rm and} \quad
(f\otimes g)\otimes h=f\otimes(g\otimes h),
\]
and an object $I$, which is a left and right unit for $\otimes$,
\[
A\otimes I=A=I\otimes A\quad {\rm and} \quad f\otimes
\mj_I=f=\mj_I\otimes f.
\]
A \emph{strict monoidal functor} between strict monoidal
categories is a functor that preserves strict monoidal structure
``on the nose'', i.e.,\ $F(A\otimes B)=F(A)\otimes F(B)$,
$F(I)=I$, etc.

Algebraist's \emph{simplicial category} $\Delta$ is an example of
strict monoidal category. The objects of $\Delta$ are all finite
ordinals $0=\emptyset$, $1=\{0\}$, $\ldots,n=\{0,\ldots,n-1\}$,
etc.  The arrows of $\Delta$ from $n$ to $m$ are all order
preserving functions from the set $n$ to the set $m$, i.e.,\
$f:n\str m$ satisfying: if $i< j$ and $i,j\in n$, then $f(i)\leq
f(j)$. We use the standard graphical presentation for arrows of
$\Delta$. For example, the unique arrows from 2 to 1 and from 0 to
1 are graphically presented as:
\begin{center}
\begin{picture}(150,40)(0,0)

\put(10,10){\circle*{2}} \put(0,30){\circle*{2}}
\put(20,30){\circle*{2}}

\put(10,2){\makebox(0,0)[b]{\scriptsize $0$}}
\put(0,33){\makebox(0,0)[b]{\scriptsize $0$}}
\put(20,33){\makebox(0,0)[b]{\scriptsize $1$}}

\put(10,10){\line(-1,2){10}} \put(10,10){\line(1,2){10}}

\put(-20,20){\makebox(0,0)[r]{$2\str 1$}}

\put(150,10){\circle*{2}} \put(150,2){\makebox(0,0)[b]{\scriptsize
$0$}}

\put(120,20){\makebox(0,0)[r]{$0\str 1$}}

\end{picture}
\end{center}

A bifunctor $\otimes:\Delta\times\Delta\str\Delta$ is defined on
objects as addition and on arrows as placing ``side by side'',
i.e.,\ for $f:n\str m$ and $f':n'\str m'$
\[
(f\otimes f')(i)=\left\{
\begin{array}{cl}
f(i), & \mbox{when}\hspace{1em} 0\leq i\leq n-1,
\\[1ex]
m+f'(i-n), & \mbox{when}\hspace{1em} n\leq i\leq n+n'-1,

\end{array}
\right .
\]
and 0 serves as the unit $I$.

A \emph{monoid} in a strict monoidal category $({\cal
M},\otimes,I)$ is a triple $(M,\mu:M\otimes M\str M,\eta:I\str M)$
such that
\[
\mu\circ(\mu\otimes\mj_M)=\mu\circ(\mj_M\otimes\mu)\quad {\rm and}
\quad \mu\circ(\mj_M\otimes\eta)=\mj_M=\mu\circ(\eta\otimes\mj_M).
\]
For example, $(1,\m1,\e1)$ is a monoid in $\Delta$, where \m1 and
\e1 are the above graphical presentations of the arrows of
$\Delta$ from 2 to 1, and from 0 to 1. The following result, taken
over from \cite[VII.5, Proposition~1]{ML71}, shows the
``universal'' property of this monoid.

{\prop Given a monoid $(M,\mu,\eta)$ in a strict monoidal category
$\cal M$, there is a unique strict monoidal functor
$F:\Delta\str{\cal M}$ such that $F(1)=M$, $F(\m1)=\mu$ and
$F(\e1)=\eta$.}

\vspace{2ex}

Let $\Delta_{par}$ be the category with the same objects as
$\Delta$, whose arrows are order preserving partial functions.
Then $(1,\m1,\e1)$ is a monoid in the strict monoidal category
$\Delta_{par}$ with the same tensor and unit as $\Delta$. The
empty partial function from 1 to 0 is graphically presented as
\p1. By \cite[Proposition~6.2]{PT13} we have the following
universal property of this monoid.

{\prop Given a monoid $(M,\mu,\eta)$ in a strict monoidal category
$\cal M$ whose monoidal structure is given by finite products,
there is a unique strict monoidal functor $F:\Delta_{par}\str{\cal
M}$ such that $F(1)=M$, $F(\m1)=\mu$, $F(\e1)=\eta$ and $F(\p1)$
is the unique arrow from $M$ to the unit (a terminal object of
$\cal M$).}

\vspace{2ex}

Topologist's simplicial category is the full subcategory of
$\Delta$ on nonempty ordinals as objects. We identify this
category with the subcategory of \top. The object $n+1$ is
identified with the standard ordered simplex
\[
\Delta^n=\{(t_0,\ldots,t_n)\mid t_0,\ldots, t_n\geq 0, \sum_i
t_i=1\},
\]
and an arrow $f\!:n\pl 1\str m\pl 1$ is identified with the affine
map defined by
\[
f(t_0,\ldots,t_n)=(s_0,\ldots,s_m),\; \mbox{\rm where }
s_j=\sum_{f(i)=j}t_i.
\]

We denote by $\Delta^{op}$ the opposite of topologist's simplicial
category and rename its objects so that the ordinal $n+1$ is
denoted by $[n]$, i.e.,\ $[n]=\{0,\ldots,n\}$. Let $\Delta_{Int}$
be the subcategory of $\Delta$ whose objects are finite ordinals
greater or equal to 2 and whose arrows are interval maps, i.e.,\
order-preserving functions, which preserve, moreover, the first
and the last element.

The categories $\Delta^{op}$ and $\Delta_{Int}$ are isomorphic by
the functor $\cal J$ (see \cite[Section~6]{PT13}). This functor
maps the object $[n]$ to $n+2$ and it maps the generating arrows
of $\Delta^{op}$ in the following way.

\begin{center}
\begin{picture}(220,50)(100,0)

\put(100,10){\circle*{2}} \put(130,10){\circle*{2}}
\put(145,10){\circle*{2}} \put(170,10){\circle*{2}}
\put(100,30){\circle*{2}} \put(130,30){\circle*{2}}
\put(145,30){\circle{2}} \put(160,30){\circle*{2}}
\put(185,30){\circle*{2}}

\put(100,0){\makebox(0,0)[b]{\scriptsize $0$}}
\put(130,0){\makebox(0,0)[b]{\scriptsize $i\mn 1$}}
\put(145,0){\makebox(0,0)[b]{\scriptsize $i$}}
\put(170,0){\makebox(0,0)[b]{\scriptsize $n\mn 1$}}

\put(100,35){\makebox(0,0)[b]{\scriptsize $0$}}
\put(130,35){\makebox(0,0)[b]{\scriptsize $i\mn 1$}}
\put(145,35){\makebox(0,0)[b]{\scriptsize $i$}}
\put(160,35){\makebox(0,0)[b]{\scriptsize $i\pl 1$}}
\put(185,35){\makebox(0,0)[b]{\scriptsize $n$}}

\put(100,29){\line(0,-1){18}} \put(130,29){\line(0,-1){18}}

\put(145.7,10.7){\line(3,4){14}} \put(170.7,10.7){\line(3,4){14}}

\put(116,20){\makebox(0,0){\ldots}}
\put(166,20){\makebox(0,0){\ldots}}

\put(200,20){\makebox(0,0)[l]{$\mapsto$}}

\put(220,10){\circle*{2}} \put(250,10){\circle*{2}}
\put(265,10){\circle*{2}} \put(280,10){\circle*{2}}
\put(305,10){\circle*{2}} \put(220,30){\circle*{2}}
\put(250,30){\circle*{2}} \put(265,30){\circle*{2}}
\put(280,30){\circle*{2}} \put(295,30){\circle*{2}}
\put(320,30){\circle*{2}}

\put(220,0){\makebox(0,0)[b]{\scriptsize $0$}}
\put(265,0){\makebox(0,0)[b]{\scriptsize $i$}}
\put(305,0){\makebox(0,0)[b]{\scriptsize $n$}}

\put(220,35){\makebox(0,0)[b]{\scriptsize $0$}}
\put(265,35){\makebox(0,0)[b]{\scriptsize $i$}}
\put(280,35){\makebox(0,0)[b]{\scriptsize $i\pl 1$}}
\put(320,35){\makebox(0,0)[b]{\scriptsize $n\pl 1$}}

\put(220,29){\line(0,-1){18}} \put(250,29){\line(0,-1){18}}
\put(265,29){\line(0,-1){18}}

\put(265.7,10.7){\line(3,4){14}} \put(280.7,10.7){\line(3,4){14}}
\put(305.7,10.7){\line(3,4){14}}

\put(236,20){\makebox(0,0){\ldots}}
\put(301,20){\makebox(0,0){\ldots}}

\end{picture}
\end{center}

\vspace{1ex}

\begin{center}
\begin{picture}(220,40)(100,0)

\put(100,30){\circle*{2}} \put(130,30){\circle*{2}}
\put(145,30){\circle*{2}} \put(170,30){\circle*{2}}
\put(100,10){\circle*{2}} \put(130,10){\circle*{2}}
\put(145,10){\circle*{2}} \put(160,10){\circle*{2}}
\put(185,10){\circle*{2}}

\put(100,35){\makebox(0,0)[b]{\scriptsize $0$}}
\put(130,35){\makebox(0,0)[b]{\scriptsize $i\mn 1$}}
\put(145,35){\makebox(0,0)[b]{\scriptsize $i$}}
\put(170,35){\makebox(0,0)[b]{\scriptsize $n\mn 1$}}

\put(100,0){\makebox(0,0)[b]{\scriptsize $0$}}
\put(130,0){\makebox(0,0)[b]{\scriptsize $i\mn 1$}}
\put(145,0){\makebox(0,0)[b]{\scriptsize $i$}}
\put(160,0){\makebox(0,0)[b]{\scriptsize $i\pl 1$}}
\put(185,0){\makebox(0,0)[b]{\scriptsize $n$}}

\put(100,29){\line(0,-1){18}} \put(130,29){\line(0,-1){18}}
\put(145,29){\line(0,-1){18}}

\put(145.7,29.3){\line(3,-4){14}}
\put(170.7,29.3){\line(3,-4){14}}
\put(116,20){\makebox(0,0){\ldots}}
\put(166,20){\makebox(0,0){\ldots}}

\put(200,20){\makebox(0,0)[l]{$\mapsto$}}

\put(220,30){\circle*{2}} \put(265,30){\circle*{2}}
\put(280,30){\circle*{2}} \put(305,30){\circle*{2}}
\put(220,10){\circle*{2}} \put(265,10){\circle*{2}}
\put(280,10){\circle{2}} \put(295,10){\circle*{2}}
\put(320,10){\circle*{2}}

\put(220,35){\makebox(0,0)[b]{\scriptsize $0$}}
\put(265,35){\makebox(0,0)[b]{\scriptsize $i$}}
\put(305,35){\makebox(0,0)[b]{\scriptsize $n$}}

\put(220,0){\makebox(0,0)[b]{\scriptsize $0$}}
\put(265,0){\makebox(0,0)[b]{\scriptsize $i$}}
\put(280,0){\makebox(0,0)[b]{\scriptsize $i\pl 1$}}
\put(320,0){\makebox(0,0)[b]{\scriptsize $n\pl 1$}}

\put(220,29){\line(0,-1){18}} \put(265,29){\line(0,-1){18}}

\put(280.7,29.3){\line(3,-4){14}}
\put(305.7,29.3){\line(3,-4){14}}

\put(244,20){\makebox(0,0){\ldots}}
\put(301,20){\makebox(0,0){\ldots}}

\end{picture}
\end{center}
The functor $\cal J$ may be visualized as the following embedding
of $\Delta^{op}$ into $\Delta$. (I am grateful to Matija Ba\v si\'
c for this remark.)

\begin{center}
\begin{picture}(230,35)(0,15)

\put(0,30){\makebox(0,0){$\textcolor{red}{\Delta^{op}}\hookrightarrow\Delta$}}

\put(50,30){\makebox(0,0){\ldots}}

\put(70,30){\textcolor{red}{\circle*{2}}}
\put(110,30){\textcolor{red}{\circle*{2}}}
\put(150,30){\textcolor{red}{\circle*{2}}}
\put(190,30){\circle*{2}} \put(230,30){\circle*{2}}

\put(90,40){\makebox(0,0){\textcolor{red}{$\rightarrow$}}}
\put(90,30){\makebox(0,0){\textcolor{red}{$\rightarrow$}}}
\put(90,20){\makebox(0,0){\textcolor{red}{$\rightarrow$}}}

\put(90,45){\makebox(0,0){$\leftarrow$}}
\put(90,35){\makebox(0,0){\textcolor{red}{$\leftarrow$}}}
\put(90,25){\makebox(0,0){\textcolor{red}{$\leftarrow$}}}
\put(90,15){\makebox(0,0){$\leftarrow$}}

\put(130,35){\makebox(0,0){\textcolor{red}{$\rightarrow$}}}
\put(130,25){\makebox(0,0){\textcolor{red}{$\rightarrow$}}}

\put(130,40){\makebox(0,0){$\leftarrow$}}
\put(130,30){\makebox(0,0){\textcolor{red}{$\leftarrow$}}}
\put(130,20){\makebox(0,0){$\leftarrow$}}

\put(170,30){\makebox(0,0){$\rightarrow$}}

\put(170,35){\makebox(0,0){$\leftarrow$}}
\put(170,25){\makebox(0,0){$\leftarrow$}}

\put(210,30){\makebox(0,0){$\leftarrow$}}

\put(230,25){\makebox(0,0)[t]{\scriptsize $0$}}
\put(190,25){\makebox(0,0)[t]{\scriptsize $1$}}
\put(150,25){\makebox(0,0)[t]{\scriptsize $2$}}
\put(110,25){\makebox(0,0)[t]{\scriptsize $3$}}
\put(70,25){\makebox(0,0)[t]{\scriptsize $4$}}

\put(150,35){\makebox(0,0)[b]{\textcolor{red}{\scriptsize $[0]$}}}
\put(110,35){\makebox(0,0)[b]{\textcolor{red}{\scriptsize $[1]$}}}
\put(70,35){\makebox(0,0)[b]{\textcolor{red}{\scriptsize $[2]$}}}

\put(280,30){\makebox(0,0){$(1)$}}

\end{picture}
\end{center}
Throughout this paper, we represent the arrows of $\Delta^{op}$ by
the graphical presentations of their $\cal J$ images in
$\Delta_{Int}$.

We have a functor ${{\cal H}\!:\Delta_{Int}\str\Delta_{par}}$
defined on objects as ${\cal H}(n)=n\mn 2$, and on arrows, for
$f\!:n\str m$, as
\[
{\cal H}(f)=
\begin{picture}(100,25)(-10,0) \put(20,0){\circle*{2}}
\put(50,0){\circle*{2}} \put(0,20){\circle{2}}
\put(70,20){\circle{2}} \put(20,20){\circle*{2}}
\put(50,20){\circle*{2}}

\put(20,-3){\makebox(0,0)[t]{\scriptsize $1$}}
\put(50,-3){\makebox(0,0)[t]{\scriptsize $m\mn 3$}}

\put(0,23){\makebox(0,0)[b]{\scriptsize $0$}}
\put(20,23){\makebox(0,0)[b]{\scriptsize $1$}}
\put(50,23){\makebox(0,0)[b]{\scriptsize $m\mn 2$}}
\put(75,23){\makebox(0,0)[b]{\scriptsize $m\mn 1$}}

\put(20,0){\line(0,1){20}} \put(50,0){\line(0,1){20}}

\put(35,10){\makebox(0,0){\ldots}}
\end{picture}
\circ f\circ
\begin{picture}(100,25)(-10,0)
\put(20,0){\circle*{2}} \put(50,0){\circle*{2}}
\put(0,0){\circle{2}} \put(70,0){\circle{2}}
\put(20,20){\circle*{2}} \put(50,20){\circle*{2}}

\put(0,-3){\makebox(0,0)[t]{\scriptsize $0$}}
\put(20,-3){\makebox(0,0)[t]{\scriptsize $1$}}
\put(50,-3){\makebox(0,0)[t]{\scriptsize $n\mn 2$}}
\put(75,-3){\makebox(0,0)[t]{\scriptsize $n\mn 1$}}

\put(20,23){\makebox(0,0)[b]{\scriptsize $0$}}
\put(50,23){\makebox(0,0)[b]{\scriptsize $n\mn 3$}}

\put(20,0){\line(0,1){20}} \put(50,0){\line(0,1){20}}

\put(35,10){\makebox(0,0){\ldots}}
\end{picture}
\]
(Intuitively, ${\cal H}(f)$ is obtained by omitting the vertices
$0$, $n\mn 1$ from the source, and $0$, $m\mn 1$ from the target
in the graphical presentation of $f$ together with all the edges
incident to these vertices.) It is not difficult to see that
${\cal H}(\mj_n)=\mj_{n-2}$, and that for a pair of arrows
$f\!:n\str m$ and $g\!:m\str k$ of $\Delta_{Int}$ we have
\[
{\cal H}(g)\circ {\cal H}(f)(i)=\left\{
\begin{array}{ll}
g(f(i+1))-1, & f(i\pl 1)\not\in\{0,m\mn 1 \}\wedge g(f(i\pl
1))\not\in \{0,k\mn 1 \}
\\[1ex]
{\rm undefined}, & {\rm otherwise},
\end{array}
\right .
\]
and
\[
{\cal H}(g\circ f)(i)=\left\{
\begin{array}{ll}
g(f(i+1))-1, & g(f(i\pl 1))\not\in \{0,k\mn 1 \}
\\[1ex]
{\rm undefined}, & {\rm otherwise}.
\end{array}
\right .
\]
Since $g(f(i\pl 1))\not\in \{0,k\mn 1 \}$ implies $f(i\pl
1)\not\in\{0,m\mn 1 \}$, we have that ${\cal H}(g)\circ {\cal
H}(f)(i)={\cal H}(g\circ f)(i)$, and ${\cal H}$ so defined is
indeed a functor.

A \emph{simplicial object} $X$ in a category $\cal M$ is a functor
$X:\Delta^{op}\str{\cal M}$. The following proposition is a
corollary of Proposition~2.2.

{\prop Given a monoid $(M,\mu,\eta)$ in a strict monoidal category
$\cal M$ whose monoidal structure is given by finite products,
there is a simplicial object $X$ in $\cal M$ such that
$X([n])=M^n$, $X(\mm2)=\mu$, $X(\ee2)=\eta$.}

\vspace{2ex}

\dkz Take $X$ to be the composition $F\circ {\cal H}\circ {\cal
J}$, for $F$ as in Proposition~2.2. \qed

\vspace{2ex}

Note that both \pp2 and \ppp3 are mapped by $X$ to the unique
arrow from $M$ to the unit $M^0$ (a terminal object of $\cal M$).
We say that a simplicial object in $\cal M$ obtained as the
composition $F\circ {\cal H}\circ {\cal J}$, for $F$ as in
Proposition~2.2, is the \emph{reduced bar construction based} on
$M$ (see \cite{T79} and \cite{PT13}).

For $X$ a simplicial object, we abbreviate $X([n])$ by $X_n$.
Also, for $f$ an arrow of $\Delta^{op}$, we abbreviate $X(f)$ by
$f$ whenever the simplicial object $X$ is determined by the
context.

For $n\geq 2$, consider the arrows $i_1,\ldots, i_n:[n]\str [1]$
of $\Delta^{op}$ graphically presented as follows.
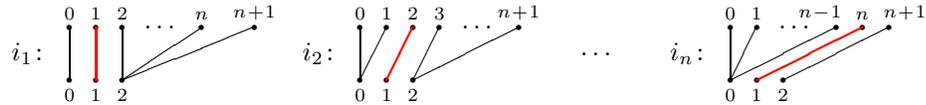
\begin{center}
\begin{picture}(320,40)(40,0)

\put(40,20){\makebox(0,0)[r]{$i_1\!:$}}

\put(50,10){\circle*{2}} \put(60,10){\circle*{2}}
\put(70,10){\circle*{2}} \put(50,30){\circle*{2}}
\put(60,30){\circle*{2}} \put(70,30){\circle*{2}}
\put(85,30){\makebox(0,0){\ldots}} \put(100,30){\circle*{2}}
\put(120,30){\circle*{2}}

\put(50,2){\makebox(0,0)[b]{\scriptsize $0$}}
\put(60,2){\makebox(0,0)[b]{\scriptsize $1$}}
\put(70,2){\makebox(0,0)[b]{\scriptsize $2$}}
\put(50,33){\makebox(0,0)[b]{\scriptsize $0$}}
\put(60,33){\makebox(0,0)[b]{\scriptsize $1$}}
\put(70,33){\makebox(0,0)[b]{\scriptsize $2$}}
\put(100,33){\makebox(0,0)[b]{\scriptsize $n$}}
\put(120,33){\makebox(0,0)[b]{\scriptsize $n\!+\!1$}}

{\thinlines \put(50,10){\line(0,1){20}}
\put(70,10){\line(0,1){20}} \put(70,10){\line(3,2){30}}
\put(70,10){\line(5,2){50}}} {\thicklines
\put(60,10){\textcolor{red}{\line(0,1){20}}}}

\put(150,20){\makebox(0,0)[r]{$i_2\!:$}}

\put(160,10){\circle*{2}} \put(170,10){\circle*{2}}
\put(180,10){\circle*{2}} \put(160,30){\circle*{2}}
\put(170,30){\circle*{2}} \put(180,30){\circle*{2}}
\put(190,30){\circle*{2}} \put(205,30){\makebox(0,0){\ldots}}
\put(220,30){\circle*{2}}

\put(160,2){\makebox(0,0)[b]{\scriptsize $0$}}
\put(170,2){\makebox(0,0)[b]{\scriptsize $1$}}
\put(180,2){\makebox(0,0)[b]{\scriptsize $2$}}
\put(160,33){\makebox(0,0)[b]{\scriptsize $0$}}
\put(170,33){\makebox(0,0)[b]{\scriptsize $1$}}
\put(180,33){\makebox(0,0)[b]{\scriptsize $2$}}
\put(190,33){\makebox(0,0)[b]{\scriptsize $3$}}
\put(220,33){\makebox(0,0)[b]{\scriptsize $n\!+\!1$}}

{\thinlines \put(160,10){\line(0,1){20}}
\put(160,10){\line(1,2){10}} \put(180,10){\line(1,2){10}}
\put(180,10){\line(2,1){40}}} {\thicklines
\put(170,10){\textcolor{red}{\line(1,2){10}}}}

\put(250,20){\makebox(0,0){\ldots}}

\put(290,20){\makebox(0,0)[r]{$i_n\!:$}}

\put(300,10){\circle*{2}} \put(310,10){\circle*{2}}
\put(320,10){\circle*{2}} \put(300,30){\circle*{2}}
\put(310,30){\circle*{2}} \put(340,30){\circle*{2}}
\put(350,30){\circle*{2}} \put(325,30){\makebox(0,0){\ldots}}
\put(360,30){\circle*{2}}

\put(300,2){\makebox(0,0)[b]{\scriptsize $0$}}
\put(310,2){\makebox(0,0)[b]{\scriptsize $1$}}
\put(320,2){\makebox(0,0)[b]{\scriptsize $2$}}
\put(300,33){\makebox(0,0)[b]{\scriptsize $0$}}
\put(310,33){\makebox(0,0)[b]{\scriptsize $1$}}
\put(342,33){\makebox(0,0)[br]{\scriptsize $n\!-\!1$}}
\put(350,33){\makebox(0,0)[b]{\scriptsize $n$}}
\put(358,33){\makebox(0,0)[bl]{\scriptsize $n\!+\!1$}}

{\thinlines \put(300,10){\line(0,1){20}}
\put(300,10){\line(1,2){10}} \put(300,10){\line(2,1){40}}
\put(320,10){\line(2,1){40}}} {\thicklines
\put(310,10){\textcolor{red}{\line(2,1){40}}}}

\end{picture}
\end{center}
(It would be more appropriate to denote these arrows by
$i^n_1,\ldots, i^n_n$, but we omit the upper indices taking them
for granted.)

For arrows $f\!:C\str A$ and $g\!:C\str B$ of a strict monoidal
category $\cal M$ whose monoidal structure is given by finite
products, we denote by $\langle f,g\rangle:C\str A\times B$ the
arrow obtained by the universal property of product in $\cal M$.
For $X$ a simplicial object in $\cal M$, we denote by $p_0$ the
unique arrow from $X_0$ to the unit, i.e.,\ a terminal object
$(X_1)^0$ of $\cal M$, and we denote by $p_1$ the identity arrow
from $X_1$ to $X_1$. For $n\geq 2$ and the above mentioned arrows
$i_1,\ldots, i_n:[n]\str [1]$ of $\Delta^{op}$, we denote by $p_n$
the arrow
\[
\langle i_1,\ldots,i_n\rangle\!:X_n\str (X_1)^n,
\]
where by our convention, $i_j$ abbreviates $X(i_j)$.

Let $X$ be the reduced bar construction based on a monoid $M$.
Since $X_0$ is the unit $M^0$ and for $n\geq 2$, the arrow
$i_j\!:M^n\str M$ is the $j$th projection, we have that for every
$n\geq 0$, the arrow $p_n$ is the identity. We show that this
property characterizes the reduced bar construction based on a
monoid in $\cal M$.

{\prop Let $\cal M$ be a strict monoidal category whose monoidal
structure is given by finite products. A simplicial object $X$ in
$\cal M$ is the reduced bar construction based on a monoid
in~$\cal M$ if and only if for every $n\geq 0$, the arrow
$p_n\!:X_n\!\str\! (X_1)^n$ is the identity.}

\vspace{2ex}

\dkz The ``only if'' part of the proof is given in the paragraph
preceding this proposition. For the ``if'' part of the proof, let
us denote $X_1$ by $M$. By our convention, the $X$ images of
arrows of $\Delta^{op}$ are denoted just by their names or
graphical presentations. We show that
\[
(M,\mm2,\ee2)
\]
is a monoid in~$\cal M$. Let $k^1_{M^2,M}\!:M^2\times M\str M^2$
and $k^2_{M^2,M}\!:M^2\times M\str M$ be the first and the second
projection respectively. Since $p_3=\langle
i_1,i_2,i_3\rangle\!:M^3\str M^3$ is the identity, we have that
$k^1_{M^2,M}=\langle i_1,i_2\rangle$ and $k^2_{M^2,M}=i_3=\i3$.

For arrows $f\!:C\str A$, $g\!:C\str B$, $h\!:D\str C$,
$f_1\!:A_1\str B_1$, $f_2\!:A_2\str B_2$ and projections
$k^1_{A_1,A_2}\!:A_1\times A_2\str A_1$ and
$k^2_{A_1,A_2}\!:A_1\times A_2\str A_2$, the following equations
hold in $\cal M$
\[
\langle f\circ h,g\circ h \rangle=\langle f,g\rangle\circ
h,\quad\quad f_1\times f_2=\langle f_1\circ k^1_{A_1,A_2},f_2\circ
k^2_{A_1,A_2}\rangle.
\]

We have
\[
k^1_{M^2,M}=\langle i_1,i_2\rangle=
\langle\begin{picture}(34,10)(0,2) \put(23,0){\line(-1,2){5}}
\put(23,0){\line(1,2){5}} \put(23,0){\line(0,1){10}}
\put(2,0){\line(0,1){10}} \put(10,0){\line(0,1){10}}
\put(23,0){\circle*{1.5}} \put(23,10){\circle*{1.5}}
\put(18,10){\circle*{1.5}} \put(28,10){\circle*{1.5}}
\put(2,0){\circle*{1.5}} \put(2,10){\circle*{1.5}}
\put(10,0){\circle*{1.5}} \put(10,10){\circle*{1.5}}\end{picture},
\begin{picture}(40,10)(0,2) \put(7,0){\line(-1,2){5}}
\put(7,0){\line(1,2){5}} \put(20,0){\line(0,1){10}}
\put(33,0){\line(-1,2){5}} \put(33,0){\line(1,2){5}}
\put(7,0){\circle*{1.5}} \put(2,10){\circle*{1.5}}
\put(12,10){\circle*{1.5}} \put(28,10){\circle*{1.5}}
\put(20,0){\circle*{1.5}} \put(20,10){\circle*{1.5}}
\put(33,0){\circle*{1.5}}
\put(38,10){\circle*{1.5}}\end{picture}\rangle=
\langle \begin{picture}(38,10)(0,7) \put(2,0){\line(0,1){10}}
\put(12,0){\line(0,1){10}} \put(25,0){\line(-1,2){5}}
\put(25,0){\line(1,2){5}} \put(2,10){\line(0,1){10}}
\put(12,10){\line(0,1){10}} \put(20,10){\line(0,1){10}}
\put(30,10){\line(-1,2){5}} \put(30,10){\line(1,2){5}}
\put(2,0){\circle*{1.5}} \put(2,10){\circle*{1.5}}
\put(12,0){\circle*{1.5}} \put(12,10){\circle*{1.5}}
\put(25,0){\circle*{1.5}} \put(20,10){\circle*{1.5}}
\put(30,10){\circle*{1.5}} \put(2,20){\circle*{1.5}}
\put(12,20){\circle*{1.5}} \put(20,20){\circle*{1.5}}
\put(25,20){\circle*{1.5}}
\put(35,20){\circle*{1.5}}\end{picture},
\begin{picture}(38,10)(0,7) \put(7,0){\line(-1,2){5}}
\put(7,0){\line(1,2){5}} \put(20,0){\line(0,1){10}}
\put(30,0){\line(0,1){10}} \put(2,10){\line(0,1){10}}
\put(12,10){\line(0,1){10}} \put(20,10){\line(0,1){10}}
\put(30,10){\line(-1,2){5}} \put(30,10){\line(1,2){5}}
\put(7,0){\circle*{1.5}} \put(2,10){\circle*{1.5}}
\put(20,0){\circle*{1.5}} \put(12,10){\circle*{1.5}}
\put(30,0){\circle*{1.5}} \put(20,10){\circle*{1.5}}
\put(30,10){\circle*{1.5}} \put(2,20){\circle*{1.5}}
\put(12,20){\circle*{1.5}} \put(20,20){\circle*{1.5}}
\put(25,20){\circle*{1.5}}
\put(35,20){\circle*{1.5}}\end{picture}\rangle=p_2 \circ \k1.
\]
Hence, $k^1_{M^2,M}=\k1$. Analogously, we prove that
$k^2_{M,M^2}=\kk2$. Also,
\[
\mu\times\mj=\langle \mu\circ k^1_{M^2,M},k^2_{M^2,M}\rangle =
\langle
\begin{picture}(38,10)(-2,7) \put(0,0){\line(0,1){10}}
\put(30,0){\line(0,1){10}} \put(15,0){\line(-1,2){5}}
\put(15,0){\line(1,2){5}} \put(0,10){\line(0,1){10}}
\put(10,10){\line(0,1){10}} \put(20,10){\line(0,1){10}}
\put(30,10){\line(-1,2){5}} \put(30,10){\line(1,2){5}}
\put(0,0){\circle*{1.5}} \put(0,10){\circle*{1.5}}
\put(30,0){\circle*{1.5}} \put(30,10){\circle*{1.5}}
\put(15,0){\circle*{1.5}} \put(20,10){\circle*{1.5}}
\put(10,10){\circle*{1.5}} \put(0,20){\circle*{1.5}}
\put(10,20){\circle*{1.5}} \put(20,20){\circle*{1.5}}
\put(25,20){\circle*{1.5}}
\put(35,20){\circle*{1.5}}\end{picture}, \i3\rangle = \langle
\begin{picture}(34,10)(-2,7) \put(0,0){\line(0,1){10}}
\put(10,0){\line(0,1){10}} \put(25,0){\line(-1,2){5}}
\put(25,0){\line(1,2){5}} \put(0,10){\line(0,1){10}}
\put(20,10){\line(0,1){10}} \put(30,10){\line(0,1){10}}
\put(10,10){\line(-1,2){5}} \put(10,10){\line(1,2){5}}
\put(0,0){\circle*{1.5}} \put(0,10){\circle*{1.5}}
\put(10,0){\circle*{1.5}} \put(10,10){\circle*{1.5}}
\put(25,0){\circle*{1.5}} \put(20,10){\circle*{1.5}}
\put(30,10){\circle*{1.5}} \put(0,20){\circle*{1.5}}
\put(5,20){\circle*{1.5}} \put(20,20){\circle*{1.5}}
\put(15,20){\circle*{1.5}}
\put(30,20){\circle*{1.5}}\end{picture},
\begin{picture}(34,10)(-2,7) \put(20,0){\line(0,1){10}}
\put(30,0){\line(0,1){10}} \put(5,0){\line(-1,2){5}}
\put(5,0){\line(1,2){5}} \put(0,10){\line(0,1){10}}
\put(20,10){\line(0,1){10}} \put(30,10){\line(0,1){10}}
\put(10,10){\line(-1,2){5}} \put(10,10){\line(1,2){5}}
\put(5,0){\circle*{1.5}} \put(0,10){\circle*{1.5}}
\put(20,0){\circle*{1.5}} \put(10,10){\circle*{1.5}}
\put(30,0){\circle*{1.5}} \put(20,10){\circle*{1.5}}
\put(30,10){\circle*{1.5}} \put(0,20){\circle*{1.5}}
\put(5,20){\circle*{1.5}} \put(20,20){\circle*{1.5}}
\put(15,20){\circle*{1.5}} \put(30,20){\circle*{1.5}}\end{picture}
\rangle = p_2\circ \begin{picture}(34,10)(-2,2)
\put(0,0){\line(0,1){10}} \put(20,0){\line(0,1){10}}
\put(30,0){\line(0,1){10}} \put(10,0){\line(-1,2){5}}
\put(10,0){\line(1,2){5}} \put(0,0){\circle*{1.5}}
\put(10,0){\circle*{1.5}} \put(20,0){\circle*{1.5}}
\put(30,0){\circle*{1.5}} \put(0,10){\circle*{1.5}}
\put(5,10){\circle*{1.5}} \put(20,10){\circle*{1.5}}
\put(15,10){\circle*{1.5}}
\put(30,10){\circle*{1.5}}\end{picture}.
\]
Hence, $\mu\times\mj= \begin{picture}(34,10)(-2,2)
\put(0,0){\line(0,1){10}} \put(20,0){\line(0,1){10}}
\put(30,0){\line(0,1){10}} \put(10,0){\line(-1,2){5}}
\put(10,0){\line(1,2){5}} \put(0,0){\circle*{1.5}}
\put(10,0){\circle*{1.5}} \put(20,0){\circle*{1.5}}
\put(30,0){\circle*{1.5}} \put(0,10){\circle*{1.5}}
\put(5,10){\circle*{1.5}} \put(20,10){\circle*{1.5}}
\put(15,10){\circle*{1.5}}
\put(30,10){\circle*{1.5}}\end{picture}$. Analogously, we prove
that $\mj\times \mu= \begin{picture}(34,10)(-2,2)
\put(0,0){\line(0,1){10}} \put(10,0){\line(0,1){10}}
\put(30,0){\line(0,1){10}} \put(20,0){\line(-1,2){5}}
\put(20,0){\line(1,2){5}} \put(0,0){\circle*{1.5}}
\put(10,0){\circle*{1.5}} \put(20,0){\circle*{1.5}}
\put(30,0){\circle*{1.5}} \put(0,10){\circle*{1.5}}
\put(10,10){\circle*{1.5}} \put(15,10){\circle*{1.5}}
\put(25,10){\circle*{1.5}}
\put(30,10){\circle*{1.5}}\end{picture}$. Now, $\mu\circ
(\mu\times\mj)=\mu\circ (\mj\times \mu)$, since
\[
\begin{picture}(34,15)(-2,7) \put(0,0){\line(0,1){10}}
\put(30,0){\line(0,1){10}} \put(15,0){\line(-1,2){5}}
\put(15,0){\line(1,2){5}} \put(0,10){\line(0,1){10}}
\put(20,10){\line(0,1){10}} \put(30,10){\line(0,1){10}}
\put(10,10){\line(-1,2){5}} \put(10,10){\line(1,2){5}}
\put(0,0){\circle*{1.5}} \put(15,0){\circle*{1.5}}
\put(30,0){\circle*{1.5}} \put(0,10){\circle*{1.5}}
\put(10,10){\circle*{1.5}} \put(20,10){\circle*{1.5}}
\put(30,10){\circle*{1.5}} \put(0,20){\circle*{1.5}}
\put(5,20){\circle*{1.5}} \put(20,20){\circle*{1.5}}
\put(15,20){\circle*{1.5}} \put(30,20){\circle*{1.5}}\end{picture}
= \begin{picture}(34,15)(-2,7) \put(0,0){\line(0,1){10}}
\put(30,0){\line(0,1){10}} \put(15,0){\line(-1,2){5}}
\put(15,0){\line(1,2){5}} \put(0,10){\line(0,1){10}}
\put(10,10){\line(0,1){10}} \put(30,10){\line(0,1){10}}
\put(20,10){\line(-1,2){5}} \put(20,10){\line(1,2){5}}
\put(0,0){\circle*{1.5}} \put(15,0){\circle*{1.5}}
\put(30,0){\circle*{1.5}} \put(0,10){\circle*{1.5}}
\put(10,10){\circle*{1.5}} \put(20,10){\circle*{1.5}}
\put(30,10){\circle*{1.5}} \put(0,20){\circle*{1.5}}
\put(15,20){\circle*{1.5}} \put(10,20){\circle*{1.5}}
\put(25,20){\circle*{1.5}} \put(30,20){\circle*{1.5}}\end{picture}
\]

That $k^1_{M,M^0}=\mj= \begin{picture}(20,10)(0,2)
\put(2,0){\line(0,1){10}} \put(10,0){\line(0,1){10}}
\put(18,0){\line(0,1){10}} \put(2,0){\circle*{1.5}}
\put(10,0){\circle*{1.5}} \put(18,0){\circle*{1.5}}
\put(2,10){\circle*{1.5}} \put(10,10){\circle*{1.5}}
\put(18,10){\circle*{1.5}}\end{picture}$, and $k^2_{M,M^0}=\pp2$
follows from the fact that $M^0$ is the strict unit and a terminal
object of $\cal M$. Hence,
\[
\mj\times \eta=\langle k^1_{M,M^0},\eta\circ k^2_{M,M^0}\rangle=
\langle
\begin{picture}(20,10)(0,2) \put(2,0){\line(0,1){10}}
\put(10,0){\line(0,1){10}} \put(18,0){\line(0,1){10}}
\put(2,0){\circle*{1.5}} \put(10,0){\circle*{1.5}}
\put(18,0){\circle*{1.5}} \put(2,10){\circle*{1.5}}
\put(10,10){\circle*{1.5}} \put(18,10){\circle*{1.5}}
\end{picture}, \begin{picture}(20,10)(0,7) \put(7,0){\line(0,1){10}}
\put(18,0){\line(0,1){10}} \put(7,10){\line(-1,2){5}}
\put(7,10){\line(1,2){5}} \put(18,10){\line(0,1){10}}
\put(7,0){\circle*{1.5}} \put(12.5,0){\circle*{1.5}}
\put(18,0){\circle*{1.5}} \put(7,10){\circle*{1.5}}
\put(18,10){\circle*{1.5}} \put(2,20){\circle*{1.5}}
\put(12,20){\circle*{1.5}} \put(18,20){\circle*{1.5}}
\end{picture}\rangle=
\langle
\begin{picture}(30,10)(0,7) \put(2,0){\line(0,1){10}}
\put(12,0){\line(0,1){10}} \put(23,0){\line(-1,2){5}}
\put(23,0){\line(1,2){5}} \put(2,10){\line(0,1){10}}
\put(12,10){\line(0,1){10}} \put(28,10){\line(0,1){10}}
\put(2,0){\circle*{1.5}} \put(12,0){\circle*{1.5}}
\put(23,0){\circle*{1.5}} \put(2,10){\circle*{1.5}}
\put(12,10){\circle*{1.5}} \put(18,10){\circle*{1.5}}
\put(28,10){\circle*{1.5}} \put(2,20){\circle*{1.5}}
\put(12,20){\circle*{1.5}} \put(28,20){\circle*{1.5}}
\end{picture}, \begin{picture}(30,10)(0,7) \put(18,0){\line(0,1){10}}
\put(28,0){\line(0,1){10}} \put(7,0){\line(-1,2){5}}
\put(7,0){\line(1,2){5}} \put(2,10){\line(0,1){10}}
\put(12,10){\line(0,1){10}} \put(28,10){\line(0,1){10}}
\put(18,0){\circle*{1.5}} \put(28,0){\circle*{1.5}}
\put(7,0){\circle*{1.5}} \put(2,10){\circle*{1.5}}
\put(12,10){\circle*{1.5}} \put(18,10){\circle*{1.5}}
\put(28,10){\circle*{1.5}} \put(2,20){\circle*{1.5}}
\put(12,20){\circle*{1.5}} \put(28,20){\circle*{1.5}}
\end{picture}\rangle = p_2\circ \begin{picture}(28,10)(0,2) \put(2,0){\line(0,1){10}}
\put(10,0){\line(0,1){10}} \put(26,0){\line(0,1){10}}
\put(2,0){\circle*{1.5}} \put(10,0){\circle*{1.5}}
\put(18,0){\circle*{1.5}} \put(26,0){\circle*{1.5}}
\put(2,10){\circle*{1.5}} \put(10,10){\circle*{1.5}}
\put(26,10){\circle*{1.5}}
\end{picture} = \begin{picture}(28,10)(0,2) \put(2,0){\line(0,1){10}}
\put(10,0){\line(0,1){10}} \put(26,0){\line(0,1){10}}
\put(2,0){\circle*{1.5}} \put(10,0){\circle*{1.5}}
\put(18,0){\circle*{1.5}} \put(26,0){\circle*{1.5}}
\put(2,10){\circle*{1.5}} \put(10,10){\circle*{1.5}}
\put(26,10){\circle*{1.5}}
\end{picture}.
\]
Now, $\mu\circ(\mj\times\eta)=\mj$, since
\[
\begin{picture}(30,10)(0,7) \put(2,0){\line(0,1){10}}
\put(28,0){\line(0,1){10}} \put(15,0){\line(-1,2){5}}
\put(15,0){\line(1,2){5}} \put(2,10){\line(0,1){10}}
\put(10,10){\line(0,1){10}} \put(28,10){\line(0,1){10}}
\put(2,0){\circle*{1.5}} \put(15,0){\circle*{1.5}}
\put(28,0){\circle*{1.5}} \put(2,10){\circle*{1.5}}
\put(10,10){\circle*{1.5}} \put(20,10){\circle*{1.5}}
\put(28,10){\circle*{1.5}} \put(2,20){\circle*{1.5}}
\put(10,20){\circle*{1.5}} \put(28,20){\circle*{1.5}}
\end{picture} = \begin{picture}(20,10)(0,2) \put(2,0){\line(0,1){10}}
\put(10,0){\line(0,1){10}} \put(18,0){\line(0,1){10}}
\put(2,0){\circle*{1.5}} \put(10,0){\circle*{1.5}}
\put(18,0){\circle*{1.5}} \put(2,10){\circle*{1.5}}
\put(10,10){\circle*{1.5}} \put(18,10){\circle*{1.5}}
\end{picture} = \mj.
\]
Analogously, we prove that $\mu\circ(\eta\times\mj)=\mj$, and
conclude that $M$ is a monoid in~$\cal M$.

Let $Y$ be the reduced bar construction based on $M$. We show that
$X=Y$. It is clear that the object parts of the functors $X$ and
$Y$ coincide. We prove that for every arrow $f\!:[m]\str[n]$ of
$\Delta^{op}$, the arrows $X(f)$ and $Y(f)$ are equal in $\cal M$.

If $n=0$, then this is trivial since $X_0$, which is equal to
$M^0$, is a terminal object of $\cal M$. If $n=1$, then $f$ has
one of the following forms
\[
\begin{picture}(30,10)(0,2)
\put(7,0){\line(-1,2){5}} \put(7,0){\line(1,2){5}}
\put(23,0){\line(-1,2){5}} \put(23,0){\line(1,2){5}}
\put(7,0){\circle*{1.5}} \put(2,10){\circle*{1.5}}
\put(12,10){\circle*{1.5}} \put(15,0){\circle*{1.5}}
\put(23,0){\circle*{1.5}} \put(18,10){\circle*{1.5}}
\put(28,10){\circle*{1.5}} \multiput(5,10)(2,0){3}{\circle*{.5}}
\multiput(21,10)(2,0){3}{\circle*{.5}}
\end{picture}\quad\quad{\rm or}\quad\quad
\begin{picture}(36,10)(0,2)
\put(7,0){\line(-1,2){5}} \put(7,0){\line(1,2){5}}
\put(18,0){\line(0,1){10}} \put(29,0){\line(-1,2){5}}
\put(29,0){\line(1,2){5}} \put(7,0){\circle*{1.5}}
\put(2,10){\circle*{1.5}} \put(12,10){\circle*{1.5}}
\put(18,0){\circle*{1.5}} \put(18,10){\circle*{1.5}}
\put(29,0){\circle*{1.5}} \put(24,10){\circle*{1.5}}
\put(34,10){\circle*{1.5}} \multiput(5,10)(2,0){3}{\circle*{.5}}
\multiput(27,10)(2,0){3}{\circle*{.5}}
\end{picture}\quad\quad{\rm or}\quad\quad
\begin{picture}(46,10)(0,2)
\put(7,0){\line(-1,2){5}} \put(7,0){\line(1,2){5}}
\put(23,0){\line(-1,2){5}} \put(23,0){\line(1,2){5}}
\put(39,0){\line(-1,2){5}} \put(39,0){\line(1,2){5}}
\put(7,0){\circle*{1.5}} \put(2,10){\circle*{1.5}}
\put(12,10){\circle*{1.5}} \put(23,0){\circle*{1.5}}
\put(18,10){\circle*{1.5}} \put(28,10){\circle*{1.5}}
\put(39,0){\circle*{1.5}} \put(34,10){\circle*{1.5}}
\put(44,10){\circle*{1.5}} \multiput(5,10)(2,0){3}{\circle*{.5}}
\multiput(21,10)(2,0){3}{\circle*{.5}}
\multiput(37,10)(2,0){3}{\circle*{.5}}
\end{picture}
\]

In the first case, $f=\begin{picture}(30,20)(0,2)
\put(7,10){\line(-1,2){5}} \put(7,10){\line(1,2){5}}
\put(23,10){\line(-1,2){5}} \put(23,10){\line(1,2){5}}
\put(7,0){\line(0,1){10}} \put(23,0){\line(0,1){10}}
\put(7,10){\circle*{1.5}} \put(2,20){\circle*{1.5}}
\put(12,20){\circle*{1.5}} \put(23,10){\circle*{1.5}}
\put(18,20){\circle*{1.5}} \put(28,20){\circle*{1.5}}
\put(7,0){\circle*{1.5}} \put(15,0){\circle*{1.5}}
\put(23,0){\circle*{1.5}} \multiput(5,20)(2,0){3}{\circle*{.5}}
\multiput(21,20)(2,0){3}{\circle*{.5}}
\end{picture}$ and the $X$ and $Y$ images of the upper part are
equal as in the case $n=0$, while $X(\ee2)=Y(\ee2)$ by the
definition of $Y$.

In the second case, $f$ is either identity and $X(f)=Y(f)$ holds,
or $f$ is $i_j$ for some $1\leq j\leq m$. From $\langle
X(i_1),\ldots,X(i_m)\rangle=\mj$, we conclude that $X(i_j)$ is the
$j$th projection from $M^m$ to $M$. On the other hand, by the
definition of $Y$, we have that $Y(i_j)$ is the $j$th projection
from $M^m$ to $M$. Hence $X(f)=Y(f)$.

In the third case, when $f$ is \begin{picture}(56,25)(0,2)
\put(7,0){\line(-1,2){5}} \put(7,0){\line(1,2){5}}
\put(28,0){\line(-1,1){10}} \put(28,0){\line(1,1){10}}
\put(49,0){\line(-1,2){5}} \put(49,0){\line(1,2){5}}
\put(7,0){\circle*{1.5}} \put(2,10){\circle*{1.5}}
\put(12,10){\circle*{1.5}} \put(28,0){\circle*{1.5}}
\put(18,10){\circle*{1.5}} \put(38,10){\circle*{1.5}}
\put(49,0){\circle*{1.5}} \put(54,10){\circle*{1.5}}
\put(44,10){\circle*{1.5}} \multiput(5,10)(2,0){3}{\circle*{.5}}
\multiput(26,10)(2,0){3}{\circle*{.5}}
\multiput(47,10)(2,0){3}{\circle*{.5}}
\put(28,20){\makebox(0,0){$\displaystyle\overbrace{}^l$}}
\end{picture}, we proceed by induction on $l\geq 2$. In the proof
we use the fact that two arrows $g,h\!:C\str M^2$ are equal in
$\cal M$ iff $k^1_{M,M}\circ g=k^1_{M,M}\circ h$ and
$k^2_{M,M}\circ g=k^2_{M,M}\circ h$, where $k^1_{M,M}$ and
$k^2_{M,M}$ are the first and the second projection from $M^2$ to
$M$. Also, we know from above that
\[
k^1_{M,M}=X(\begin{picture}(27,10)(0,2) \put(2,0){\line(0,1){10}}
\put(10,0){\line(0,1){10}} \put(20,0){\line(-1,2){5}}
\put(20,0){\line(1,2){5}} \put(2,0){\circle*{1.5}}
\put(10,0){\circle*{1.5}} \put(20,0){\circle*{1.5}}
\put(2,10){\circle*{1.5}} \put(10,10){\circle*{1.5}}
\put(15,10){\circle*{1.5}} \put(25,10){\circle*{1.5}}
\end{picture})=Y(\begin{picture}(27,10)(0,2) \put(2,0){\line(0,1){10}}
\put(10,0){\line(0,1){10}} \put(20,0){\line(-1,2){5}}
\put(20,0){\line(1,2){5}} \put(2,0){\circle*{1.5}}
\put(10,0){\circle*{1.5}} \put(20,0){\circle*{1.5}}
\put(2,10){\circle*{1.5}} \put(10,10){\circle*{1.5}}
\put(15,10){\circle*{1.5}} \put(25,10){\circle*{1.5}}
\end{picture}),\quad
k^2_{M,M}=X(\begin{picture}(27,10)(0,2) \put(17,0){\line(0,1){10}}
\put(25,0){\line(0,1){10}} \put(7,0){\line(-1,2){5}}
\put(7,0){\line(1,2){5}} \put(17,0){\circle*{1.5}}
\put(25,0){\circle*{1.5}} \put(7,0){\circle*{1.5}}
\put(2,10){\circle*{1.5}} \put(12,10){\circle*{1.5}}
\put(17,10){\circle*{1.5}} \put(25,10){\circle*{1.5}}
\end{picture})=Y(\begin{picture}(27,10)(0,2) \put(17,0){\line(0,1){10}}
\put(25,0){\line(0,1){10}} \put(7,0){\line(-1,2){5}}
\put(7,0){\line(1,2){5}} \put(17,0){\circle*{1.5}}
\put(25,0){\circle*{1.5}} \put(7,0){\circle*{1.5}}
\put(2,10){\circle*{1.5}} \put(12,10){\circle*{1.5}}
\put(17,10){\circle*{1.5}} \put(25,10){\circle*{1.5}}
\end{picture}).
\]

If $l=2$, then $f$ is equal to \begin{picture}(46,15)(0,2)
\put(7,0){\line(0,1){10}} \put(7,10){\line(-1,2){5}}
\put(7,10){\line(1,2){5}} \put(23,0){\line(-1,2){5}}
\put(23,0){\line(1,2){5}} \put(18,10){\line(0,1){10}}
\put(28,10){\line(0,1){10}} \put(39,0){\line(0,1){10}}
\put(39,10){\line(-1,2){5}} \put(39,10){\line(1,2){5}}
\put(7,0){\circle*{1.5}} \put(7,10){\circle*{1.5}}
\put(2,20){\circle*{1.5}} \put(12,20){\circle*{1.5}}
\put(23,0){\circle*{1.5}} \put(18,10){\circle*{1.5}}
\put(28,10){\circle*{1.5}} \put(18,20){\circle*{1.5}}
\put(28,20){\circle*{1.5}} \put(39,0){\circle*{1.5}}
\put(39,10){\circle*{1.5}} \put(34,20){\circle*{1.5}}
\put(44,20){\circle*{1.5}} \multiput(5,20)(2,0){3}{\circle*{.5}}
\multiput(37,20)(2,0){3}{\circle*{.5}}
\end{picture}. Since $X(\mm2)=Y(\mm2)$, in order to prove that
$X(f)=Y(f)$, it suffices to prove that
$g=X(\begin{picture}(46,12)(0,2) \put(7,0){\line(-1,2){5}}
\put(7,0){\line(1,2){5}} \put(18,0){\line(0,1){10}}
\put(28,0){\line(0,1){10}} \put(39,0){\line(-1,2){5}}
\put(39,0){\line(1,2){5}} \put(7,0){\circle*{1.5}}
\put(2,10){\circle*{1.5}} \put(12,10){\circle*{1.5}}
\put(18,0){\circle*{1.5}} \put(28,0){\circle*{1.5}}
\put(18,10){\circle*{1.5}} \put(28,10){\circle*{1.5}}
\put(39,0){\circle*{1.5}} \put(34,10){\circle*{1.5}}
\put(44,10){\circle*{1.5}} \multiput(5,10)(2,0){3}{\circle*{.5}}
\multiput(37,10)(2,0){3}{\circle*{.5}}
\end{picture})$ is equal to $h=Y(\begin{picture}(46,12)(0,2) \put(7,0){\line(-1,2){5}}
\put(7,0){\line(1,2){5}} \put(18,0){\line(0,1){10}}
\put(28,0){\line(0,1){10}} \put(39,0){\line(-1,2){5}}
\put(39,0){\line(1,2){5}} \put(7,0){\circle*{1.5}}
\put(2,10){\circle*{1.5}} \put(12,10){\circle*{1.5}}
\put(18,0){\circle*{1.5}} \put(28,0){\circle*{1.5}}
\put(18,10){\circle*{1.5}} \put(28,10){\circle*{1.5}}
\put(39,0){\circle*{1.5}} \put(34,10){\circle*{1.5}}
\put(44,10){\circle*{1.5}} \multiput(5,10)(2,0){3}{\circle*{.5}}
\multiput(37,10)(2,0){3}{\circle*{.5}}
\end{picture})$. By relying on the second case for $\dagger$, we have that
\[
k^1_{M,M}\circ g=X(\begin{picture}(46,15)(0,7)
\put(7,0){\line(0,1){10}} \put(18,0){\line(0,1){10}}
\put(34,0){\line(-1,2){5}} \put(34,0){\line(1,2){5}}
\put(7,10){\line(-1,2){5}} \put(7,10){\line(1,2){5}}
\put(18,10){\line(0,1){10}} \put(29,10){\line(0,1){10}}
\put(39,10){\line(-1,2){5}} \put(39,10){\line(1,2){5}}
\put(7,0){\circle*{1.5}} \put(18,0){\circle*{1.5}}
\put(34,0){\circle*{1.5}} \put(7,10){\circle*{1.5}}
\put(2,20){\circle*{1.5}} \put(12,20){\circle*{1.5}}
\put(18,10){\circle*{1.5}} \put(29,10){\circle*{1.5}}
\put(18,20){\circle*{1.5}} \put(29,20){\circle*{1.5}}
\put(39,10){\circle*{1.5}} \put(34,20){\circle*{1.5}}
\put(44,20){\circle*{1.5}} \multiput(5,20)(2,0){3}{\circle*{.5}}
\multiput(37,20)(2,0){3}{\circle*{.5}}
\end{picture}) \stackrel{\dagger}{=}
Y(\begin{picture}(46,15)(0,7) \put(7,0){\line(0,1){10}}
\put(18,0){\line(0,1){10}} \put(34,0){\line(-1,2){5}}
\put(34,0){\line(1,2){5}} \put(7,10){\line(-1,2){5}}
\put(7,10){\line(1,2){5}} \put(18,10){\line(0,1){10}}
\put(29,10){\line(0,1){10}} \put(39,10){\line(-1,2){5}}
\put(39,10){\line(1,2){5}} \put(7,0){\circle*{1.5}}
\put(18,0){\circle*{1.5}} \put(34,0){\circle*{1.5}}
\put(7,10){\circle*{1.5}} \put(2,20){\circle*{1.5}}
\put(12,20){\circle*{1.5}} \put(18,10){\circle*{1.5}}
\put(29,10){\circle*{1.5}} \put(18,20){\circle*{1.5}}
\put(29,20){\circle*{1.5}} \put(39,10){\circle*{1.5}}
\put(34,20){\circle*{1.5}} \put(44,20){\circle*{1.5}}
\multiput(5,20)(2,0){3}{\circle*{.5}}
\multiput(37,20)(2,0){3}{\circle*{.5}}
\end{picture}) = k^1_{M,M}\circ h.
\]
Analogously, we prove that $k^2_{M,M}\circ g=k^2_{M,M}\circ h$.
Hence, $g=h$.

If $l>2$, then $f$ is equal to \begin{picture}(51,20)(0,2)
\put(7,0){\line(0,1){10}} \put(7,10){\line(-1,2){5}}
\put(7,10){\line(1,2){5}} \put(28,0){\line(-1,2){5}}
\put(28,0){\line(1,2){5}} \put(23,10){\line(-1,2){5}}
\put(23,10){\line(1,2){5}} \put(33,10){\line(0,1){10}}
\put(44,0){\line(0,1){10}} \put(44,10){\line(-1,2){5}}
\put(44,10){\line(1,2){5}} \put(7,0){\circle*{1.5}}
\put(7,10){\circle*{1.5}} \put(2,20){\circle*{1.5}}
\put(12,20){\circle*{1.5}} \put(28,0){\circle*{1.5}}
\put(23,10){\circle*{1.5}} \put(33,10){\circle*{1.5}}
\put(18,20){\circle*{1.5}} \put(28,20){\circle*{1.5}}
\put(33,20){\circle*{1.5}} \put(44,0){\circle*{1.5}}
\put(44,10){\circle*{1.5}} \put(39,20){\circle*{1.5}}
\put(49,20){\circle*{1.5}} \multiput(5,20)(2,0){3}{\circle*{.5}}
\multiput(21,20)(2,0){3}{\circle*{.5}}
\multiput(42,20)(2,0){3}{\circle*{.5}}
\end{picture}, and it suffices to prove that $g=X(\begin{picture}(51,15)(0,2)
\put(7,0){\line(-1,2){5}} \put(7,0){\line(1,2){5}}
\put(23,0){\line(-1,2){5}} \put(23,0){\line(1,2){5}}
\put(33,0){\line(0,1){10}} \put(44,0){\line(-1,2){5}}
\put(44,0){\line(1,2){5}} \put(7,0){\circle*{1.5}}
\put(2,10){\circle*{1.5}} \put(12,10){\circle*{1.5}}
\put(23,0){\circle*{1.5}} \put(33,0){\circle*{1.5}}
\put(18,10){\circle*{1.5}} \put(28,10){\circle*{1.5}}
\put(33,10){\circle*{1.5}} \put(44,0){\circle*{1.5}}
\put(39,10){\circle*{1.5}} \put(49,10){\circle*{1.5}}
\multiput(5,10)(2,0){3}{\circle*{.5}}
\multiput(21,10)(2,0){3}{\circle*{.5}}
\multiput(42,10)(2,0){3}{\circle*{.5}}
\end{picture})$ is equal to $h=Y(\begin{picture}(51,15)(0,2)
\put(7,0){\line(-1,2){5}} \put(7,0){\line(1,2){5}}
\put(23,0){\line(-1,2){5}} \put(23,0){\line(1,2){5}}
\put(33,0){\line(0,1){10}} \put(44,0){\line(-1,2){5}}
\put(44,0){\line(1,2){5}} \put(7,0){\circle*{1.5}}
\put(2,10){\circle*{1.5}} \put(12,10){\circle*{1.5}}
\put(23,0){\circle*{1.5}} \put(33,0){\circle*{1.5}}
\put(18,10){\circle*{1.5}} \put(28,10){\circle*{1.5}}
\put(33,10){\circle*{1.5}} \put(44,0){\circle*{1.5}}
\put(39,10){\circle*{1.5}} \put(49,10){\circle*{1.5}}
\multiput(5,10)(2,0){3}{\circle*{.5}}
\multiput(21,10)(2,0){3}{\circle*{.5}}
\multiput(42,10)(2,0){3}{\circle*{.5}}
\end{picture})$. By relying on the induction hypothesis for $\dagger$, we have that
\[
k^1_{M,M}\circ g= X(\begin{picture}(51,20)(0,7)
\put(7,0){\line(0,1){10}} \put(7,10){\line(-1,2){5}}
\put(7,10){\line(1,2){5}} \put(39,0){\line(-1,2){5}}
\put(39,0){\line(1,2){5}} \put(23,10){\line(-1,2){5}}
\put(23,10){\line(1,2){5}} \put(34,10){\line(0,1){10}}
\put(23,0){\line(0,1){10}} \put(44,10){\line(-1,2){5}}
\put(44,10){\line(1,2){5}} \put(7,0){\circle*{1.5}}
\put(7,10){\circle*{1.5}} \put(2,20){\circle*{1.5}}
\put(12,20){\circle*{1.5}} \put(23,0){\circle*{1.5}}
\put(23,10){\circle*{1.5}} \put(34,10){\circle*{1.5}}
\put(18,20){\circle*{1.5}} \put(28,20){\circle*{1.5}}
\put(34,20){\circle*{1.5}} \put(39,0){\circle*{1.5}}
\put(44,10){\circle*{1.5}} \put(39,20){\circle*{1.5}}
\put(49,20){\circle*{1.5}} \multiput(5,20)(2,0){3}{\circle*{.5}}
\multiput(21,20)(2,0){3}{\circle*{.5}}
\multiput(42,20)(2,0){3}{\circle*{.5}}
\end{picture}) \stackrel{\dagger}{=}
Y(\begin{picture}(51,20)(0,7) \put(7,0){\line(0,1){10}}
\put(7,10){\line(-1,2){5}} \put(7,10){\line(1,2){5}}
\put(39,0){\line(-1,2){5}} \put(39,0){\line(1,2){5}}
\put(23,10){\line(-1,2){5}} \put(23,10){\line(1,2){5}}
\put(34,10){\line(0,1){10}} \put(23,0){\line(0,1){10}}
\put(44,10){\line(-1,2){5}} \put(44,10){\line(1,2){5}}
\put(7,0){\circle*{1.5}} \put(7,10){\circle*{1.5}}
\put(2,20){\circle*{1.5}} \put(12,20){\circle*{1.5}}
\put(23,0){\circle*{1.5}} \put(23,10){\circle*{1.5}}
\put(34,10){\circle*{1.5}} \put(18,20){\circle*{1.5}}
\put(28,20){\circle*{1.5}} \put(34,20){\circle*{1.5}}
\put(39,0){\circle*{1.5}} \put(44,10){\circle*{1.5}}
\put(39,20){\circle*{1.5}} \put(49,20){\circle*{1.5}}
\multiput(5,20)(2,0){3}{\circle*{.5}}
\multiput(21,20)(2,0){3}{\circle*{.5}}
\multiput(42,20)(2,0){3}{\circle*{.5}}
\end{picture})= k^1_{M,M}\circ h.
\]
By relying on the second case for $\dagger$, we have that
\[
k^2_{M,M}\circ g= X(\begin{picture}(51,20)(0,7)
\put(34,0){\line(0,1){10}} \put(7,10){\line(-1,2){5}}
\put(7,10){\line(1,2){5}} \put(15,0){\line(-3,4){8}}
\put(15,0){\line(3,4){8}} \put(23,10){\line(-1,2){5}}
\put(23,10){\line(1,2){5}} \put(34,10){\line(0,1){10}}
\put(44,0){\line(0,1){10}} \put(44,10){\line(-1,2){5}}
\put(44,10){\line(1,2){5}} \put(15,0){\circle*{1.5}}
\put(7,10){\circle*{1.5}} \put(2,20){\circle*{1.5}}
\put(12,20){\circle*{1.5}} \put(34,0){\circle*{1.5}}
\put(44,0){\circle*{1.5}} \put(23,10){\circle*{1.5}}
\put(34,10){\circle*{1.5}} \put(18,20){\circle*{1.5}}
\put(28,20){\circle*{1.5}} \put(34,20){\circle*{1.5}}
\put(44,10){\circle*{1.5}} \put(39,20){\circle*{1.5}}
\put(49,20){\circle*{1.5}} \multiput(5,20)(2,0){3}{\circle*{.5}}
\multiput(21,20)(2,0){3}{\circle*{.5}}
\multiput(42,20)(2,0){3}{\circle*{.5}}
\end{picture}) \stackrel{\dagger}{=}
Y(\begin{picture}(51,20)(0,7) \put(34,0){\line(0,1){10}}
\put(7,10){\line(-1,2){5}} \put(7,10){\line(1,2){5}}
\put(15,0){\line(-3,4){8}} \put(15,0){\line(3,4){8}}
\put(23,10){\line(-1,2){5}} \put(23,10){\line(1,2){5}}
\put(34,10){\line(0,1){10}} \put(44,0){\line(0,1){10}}
\put(44,10){\line(-1,2){5}} \put(44,10){\line(1,2){5}}
\put(15,0){\circle*{1.5}} \put(7,10){\circle*{1.5}}
\put(2,20){\circle*{1.5}} \put(12,20){\circle*{1.5}}
\put(34,0){\circle*{1.5}} \put(44,0){\circle*{1.5}}
\put(23,10){\circle*{1.5}} \put(34,10){\circle*{1.5}}
\put(18,20){\circle*{1.5}} \put(28,20){\circle*{1.5}}
\put(34,20){\circle*{1.5}} \put(44,10){\circle*{1.5}}
\put(39,20){\circle*{1.5}} \put(49,20){\circle*{1.5}}
\multiput(5,20)(2,0){3}{\circle*{.5}}
\multiput(21,20)(2,0){3}{\circle*{.5}}
\multiput(42,20)(2,0){3}{\circle*{.5}}
\end{picture})= k^2_{M,M}\circ h.
\]
Hence, $g=h$. This concludes the case when $f$ maps $[m]$ to
$[1]$.

Suppose now that $f\!:[m]\str[n]$ is an arrow of $\Delta^{op}$ and
$n\geq 2$. As in the case when $n=1$, we conclude that for every
$1\leq j\leq n$,
\[
X(i_j\circ f)=Y(i_j\circ f).
\]
Since,
\[
\langle X(i_1),\ldots,X(i_n)\rangle = \langle
Y(i_1),\ldots,Y(i_n)\rangle =\mj_{M^n},
\]
we have that
\begin{tabbing}
\hspace{1.5em}$X(f)$ \= $=\langle X(i_1),\ldots,X(i_n)\rangle
\circ X(f)=\langle X(i_1)\circ X(f),\ldots,X(i_n)\circ
X(f)\rangle$
\\[1ex]
\> $=\langle Y(i_1)\circ Y(f),\ldots,Y(i_n)\circ
Y(f)\rangle=\langle Y(i_1),\ldots,Y(i_n)\rangle \circ Y(f)$
\\[1ex]
\> $=Y(f)$. \` $\dashv$
\end{tabbing}

\section{Segal's simplicial spaces}

Let \top\ be the category of compactly generated Hausdorff spaces.
For a simplicial object in \top, i.e.,\ a \emph{simplicial space}
$X$, a relaxed form of the condition
\[
\mbox{\rm for every}\; n,\; p_n\!:X_n\!\str\! (X_1)^n\; \mbox{\rm
is the identity},
\]
reads
\[
\mbox{\rm for every}\; n,\; p_n\!:X_n\!\str\! (X_1)^n\; \mbox{\rm
is a homotopy equivalence}.
\]

Segal, \cite{S74}, used simplicial spaces satisfying this relaxed
condition for his delooping constructions and we call them
\emph{Segal's simplicial spaces}. (Note that, for the sake of
simplicity, this notion is weaker than the one defined in
\cite{P14}.) Essentially as in the proof of Proposition~2.4, one
can show the following.

{\prop If $X\!:\Delta^{op}\str\top$ is Segal's simplicial space,
then $X_1$ is a homotopy associative \mbox{\emph{H}-space} whose
multiplication is given by the composition
\[(X_1)^2\stackrel{p_2^{-1}}{\longrightarrow}X_2\stackrel{d^2_1}{\longrightarrow}X_1,\]
where $p_2^{-1}$ is an arbitrary homotopy inverse to $p_2$, and
whose unit is $s^1_0(x_0)$, for an arbitrary $x_0\in X_0$.}

\vspace{1ex}

\noindent  (A complete proof of this proposition is given in
\cite[Appendix, Proof of Lemma~3.1]{P14}.)

\vspace{2ex}

The \emph{realization} of a simplicial space $X$, is the quotient
space
\[
|X|=\left(\coprod_n X_n\times\Delta^n\right){\Big\slash}\sim,
\]
where $\sim$ is the smallest equivalence relation on $\coprod_n
X_n\times\Delta^n$ such that for every $f\!:[n]\str[m]$
of~$\Delta$, $x\in X_m$ and $t\in\Delta^n$
\[
(f^{op}(x),t)\sim(x,f(t)).
\]

A \emph{simplicial map} is a natural transformation between
simplicial spaces. Note that the realization is \emph{functorial},
i.e.,\ it is defined also for simplicial maps. For simplicial
spaces $X$ and $Y$ the \emph{product} $X\times Y$ is defined so
that $(X\times Y)_n=X_n\times Y_n$ and $(X\times Y)(f)=X(f)\times
Y(f)$. The $n$th component of the first projection $k^1\!:X\times
Y\str X$ is the first projection $k^1_n:X_n\times Y_n\str X_n$ and
analogously for the second projection. The realization functor
preserves products of simplicial spaces (see
\cite[Theorem~14.3]{M67}, \cite[III.3, Theorem]{GZ67} and
\cite[Corollary~11.6]{M72}) in the sense that
\[
\langle|k^1|,|k^2|\rangle\!: |X\times Y| \str |X|\times |Y|
\]
is a homeomorphism.

The following two propositions stem from \cite[Proposition 1.5
(b)]{S74} and from \cite[Appendix, Theorem A4 (ii)]{M74} (see also
\cite[Lemma~2.11]{P14}).

{\prop Let $X\!:\Delta^{op}\str\top$ be Segal's simplicial space
such that for every $m$, $X_m$ is a CW-complex. If $X_1$ with
respect to the \emph{H}-space structure is grouplike, then
$X_1\simeq \Omega|X|$.}

\vspace{2ex}

{\prop Let $f\!:X\str Y$ be a simplicial map of simplicial spaces
such that for every $m$, $X_m$ and $Y_m$ are CW-complexes. If each
$f_m\!:X_m\str Y_m$ is a homotopy equivalence, then $|f|\!:|X|\str
|Y|$ is a homotopy equivalence.}

\section{Segal's bisimplicial spaces}

A \emph{bisimplicial space} is a functor $X\!:\Delta^{op}\times
\Delta^{op}\str\top$ and it may be visualized as the following
graph (see the red subgraph of (1)).
\begin{center}
\begin{picture}(210,110)(20,15)

\put(40,85){\makebox(0,0){$\downarrow$}}
\put(50,85){\makebox(0,0){$\downarrow$}}
\put(60,85){\makebox(0,0){$\downarrow$}}

\put(45,85){\makebox(0,0){$\uparrow$}}
\put(55,85){\makebox(0,0){$\uparrow$}}

\put(120,85){\makebox(0,0){$\downarrow$}}
\put(130,85){\makebox(0,0){$\downarrow$}}
\put(140,85){\makebox(0,0){$\downarrow$}}

\put(125,85){\makebox(0,0){$\uparrow$}}
\put(135,85){\makebox(0,0){$\uparrow$}}

\put(200,85){\makebox(0,0){$\downarrow$}}
\put(210,85){\makebox(0,0){$\downarrow$}}
\put(220,85){\makebox(0,0){$\downarrow$}}

\put(205,85){\makebox(0,0){$\uparrow$}}
\put(215,85){\makebox(0,0){$\uparrow$}}

\put(50,45){\makebox(0,0){$\uparrow$}}

\put(45,45){\makebox(0,0){$\downarrow$}}
\put(55,45){\makebox(0,0){$\downarrow$}}

\put(130,45){\makebox(0,0){$\uparrow$}}

\put(125,45){\makebox(0,0){$\downarrow$}}
\put(135,45){\makebox(0,0){$\downarrow$}}

\put(210,45){\makebox(0,0){$\uparrow$}}

\put(205,45){\makebox(0,0){$\downarrow$}}
\put(215,45){\makebox(0,0){$\downarrow$}}

\put(50,125){\makebox(0,0){$\vdots$}}
\put(130,125){\makebox(0,0){$\vdots$}}
\put(210,125){\makebox(0,0){$\vdots$}}

\put(20,105){\makebox(0,0){$\ldots$}}

\put(50,105){\makebox(0,0){$X_{22}$}}

\put(90,115){\makebox(0,0){$\rightarrow$}}
\put(90,105){\makebox(0,0){$\rightarrow$}}
\put(90,95){\makebox(0,0){$\rightarrow$}}

\put(90,110){\makebox(0,0){$\leftarrow$}}
\put(90,100){\makebox(0,0){$\leftarrow$}}

\put(130,105){\makebox(0,0){$X_{12}$}}

\put(170,110){\makebox(0,0){$\rightarrow$}}
\put(170,100){\makebox(0,0){$\rightarrow$}}

\put(170,105){\makebox(0,0){$\leftarrow$}}

\put(210,105){\makebox(0,0){$X_{02}$}}

\put(20,65){\makebox(0,0){$\ldots$}}

\put(50,65){\makebox(0,0){$X_{21}$}}

\put(90,75){\makebox(0,0){$\rightarrow$}}
\put(90,65){\makebox(0,0){$\rightarrow$}}
\put(90,55){\makebox(0,0){$\rightarrow$}}

\put(90,70){\makebox(0,0){$\leftarrow$}}
\put(90,60){\makebox(0,0){$\leftarrow$}}

\put(130,65){\makebox(0,0){$X_{11}$}}

\put(170,70){\makebox(0,0){$\rightarrow$}}
\put(170,60){\makebox(0,0){$\rightarrow$}}

\put(170,65){\makebox(0,0){$\leftarrow$}}

\put(210,65){\makebox(0,0){$X_{01}$}}

\put(20,25){\makebox(0,0){$\ldots$}}

\put(50,25){\makebox(0,0){$X_{20}$}}

\put(90,35){\makebox(0,0){$\rightarrow$}}
\put(90,25){\makebox(0,0){$\rightarrow$}}
\put(90,15){\makebox(0,0){$\rightarrow$}}

\put(90,30){\makebox(0,0){$\leftarrow$}}
\put(90,20){\makebox(0,0){$\leftarrow$}}

\put(130,25){\makebox(0,0){$X_{10}$}}

\put(170,30){\makebox(0,0){$\rightarrow$}}
\put(170,20){\makebox(0,0){$\rightarrow$}}

\put(170,25){\makebox(0,0){$\leftarrow$}}

\put(210,25){\makebox(0,0){$X_{00}$}}

\end{picture}
\end{center}

Let $Y_n$, for $n\geq 0$, be the realization of the $n$th column,
i.e.,\ $Y_n=|X_{n\,\underline{\mbox{\hspace{.9ex}}}}|$. Since the
realization is functorial, we obtain the simplicial space $Y$.
\begin{center}
\begin{picture}(210,30)(20,10)
\put(20,25){\makebox(0,0){$\ldots$}}

\put(50,25){\makebox(0,0){$Y_{2}$}}

\put(90,35){\makebox(0,0){$\rightarrow$}}
\put(90,25){\makebox(0,0){$\rightarrow$}}
\put(90,15){\makebox(0,0){$\rightarrow$}}

\put(90,30){\makebox(0,0){$\leftarrow$}}
\put(90,20){\makebox(0,0){$\leftarrow$}}

\put(130,25){\makebox(0,0){$Y_{1}$}}

\put(170,30){\makebox(0,0){$\rightarrow$}}
\put(170,20){\makebox(0,0){$\rightarrow$}}

\put(170,25){\makebox(0,0){$\leftarrow$}}

\put(210,25){\makebox(0,0){$Y_{0}$}}

\end{picture}
\end{center}
The \emph{realization} $|X|$ of the bisimplicial space $X$ is the
realization $|Y|$ of the simplicial space $Y$.

If the simplicial space $X_{1\,\underline{\mbox{\hspace{.9ex}}}}$
is Segal's, then, by Proposition~3.1, $X_{11}$ is a homotopy
associative H-space and this is the H-space structure we refer to
in the following proposition.

{\prop If $X\!:\Delta^{op}\times\Delta^{op}\str\top$ is a
bisimplicial space such that
$X_{1\,\underline{\mbox{\hspace{.9ex}}}}$ is Segal's, $X_{11}$
with respect to the \emph{H}-space structure is grouplike, for
every $m\geq 0$, $X_{\underline{\mbox{\hspace{.9ex}}}\,m}$ is
Segal's, and for every $n,m\geq 0$, $X_{nm}$ and $Y_n$ are
CW-complexes, then $X_{11}\simeq \Omega^2|X|$.}

\vspace{2ex}

\dkz Since $X_{1\,\underline{\mbox{\hspace{.9ex}}}}$ is Segal's
simplicial space such that for every $m$, $X_{1m}$ is a CW-complex
and $X_{11}$ with respect to the \emph{H}-space structure is
grouplike, by Proposition~3.2 we have that $X_{11}\simeq
\Omega|X_{1\,\underline{\mbox{\hspace{.9ex}}}}|=\Omega Y_1$.

For every $m$, $X_{\underline{\mbox{\hspace{.9ex}}}\,m}$ is
Segal's. Hence, for every $n$, the map $p_{nm}\!:X_{nm}\str
(X_{1m})^n$, is a homotopy equivalence. The map $p_{0m}$ is the
unique map from $X_{0m}$ to $(X_{1m})^0$, the map $p_{1m}$ is the
identity on $X_{1m}$, and for $n\geq 2$, the map $p_{nm}$ is
\[
\langle (i_1,m),\ldots,(i_n,m)\rangle\!:X_{nm}\str (X_{1m})^n.
\]
Also, for every $f\!:[m]\str [m']$ of $\Delta^{op}$ and every $n$
the following diagram commutes:
\begin{center}
\begin{picture}(140,50)(0,-5)

\put(0,30){\makebox(0,0){$X_{nm}$}}

\put(18,30){\vector(1,0){100}}

\put(140,30){\makebox(0,0){$(X_{1m})^n$}}

\put(68,33){\makebox(0,0)[b]{$p_{nm}$}}

\put(0,0){\makebox(0,0){$X_{nm'}$}}

\put(18,0){\vector(1,0){100}}

\put(140,0){\makebox(0,0){$(X_{1m'})^n$}}

\put(68,3){\makebox(0,0)[b]{$p_{nm'}$}}

\put(0,22){\vector(0,-1){15}} \put(140,22){\vector(0,-1){15}}
\put(-5,15){\makebox(0,0)[r]{$(n,f)$}}
\put(145,15){\makebox(0,0)[l]{$(1,f)^n$}}

\end{picture}
\end{center}
Hence, for every $n$, $p_{n \underline{\mbox{\hspace{.9ex}}}}$ is
a simplicial map.
\begin{center}
\begin{picture}(140,90)(0,-5)

\put(0,60){\makebox(0,0){$X_{n2}$}}

\put(20,60){\vector(1,0){100}}

\put(140,60){\makebox(0,0){$(X_{12})^n$}}

\put(70,63){\makebox(0,0)[b]{$p_{n2}$}}

\put(0,30){\makebox(0,0){$X_{n1}$}}

\put(20,30){\vector(1,0){100}}

\put(140,30){\makebox(0,0){$(X_{11})^n$}}

\put(70,33){\makebox(0,0)[b]{$p_{n1}$}}

\put(0,0){\makebox(0,0){$X_{n0}$}}

\put(20,0){\vector(1,0){100}}

\put(140,0){\makebox(0,0){$(X_{10})^n$}}

\put(70,3){\makebox(0,0)[b]{$p_{n0}$}}

\put(-5,45){\makebox(0,0){$\uparrow$}}
\put(5,45){\makebox(0,0){$\uparrow$}}
\put(-10,45){\makebox(0,0){$\downarrow$}}
\put(0,45){\makebox(0,0){$\downarrow$}}
\put(10,45){\makebox(0,0){$\downarrow$}}

\put(135,45){\makebox(0,0){$\uparrow$}}
\put(145,45){\makebox(0,0){$\uparrow$}}
\put(130,45){\makebox(0,0){$\downarrow$}}
\put(140,45){\makebox(0,0){$\downarrow$}}
\put(150,45){\makebox(0,0){$\downarrow$}}

\put(0,15){\makebox(0,0){$\uparrow$}}
\put(-5,15){\makebox(0,0){$\downarrow$}}
\put(5,15){\makebox(0,0){$\downarrow$}}

\put(140,15){\makebox(0,0){$\uparrow$}}
\put(135,15){\makebox(0,0){$\downarrow$}}
\put(145,15){\makebox(0,0){$\downarrow$}}

\put(0,78){\makebox(0,0){$\vdots$}}
\put(140,78){\makebox(0,0){$\vdots$}}

\end{picture}
\end{center}

Every $(X_{1m})^n$ is a CW-complex since the product of
CW-complexes in \top\ is a CW-complex. By Proposition~3.3, for
every $n$, $|p_{n \underline{\mbox{\hspace{.9ex}}}}|\!: Y_n\str
|(X_{1\,\underline{\mbox{\hspace{.9ex}}}})^n|$ is a homotopy
equivalence. Since $|(X_{1\,\underline{\mbox{\hspace{.9ex}}}})^0|$
is a singleton it is homeomorphic to $(Y_1)^0$ and we have that
$p_0\!:Y_0\str (Y_1)^0$, as a composition of a homeomorphism with
$|p_{0 \underline{\mbox{\hspace{.9ex}}}}|$, is a homotopy
equivalence. The map $p_1\!:Y_1\str Y_1$ is the identity. For
$n\geq 2$, $\langle
|k^1|,\ldots,|k^n|\rangle\!:|(X_{1\,\underline{\mbox{\hspace{.9ex}}}})^n|\str
|X_{1\,\underline{\mbox{\hspace{.9ex}}}}|^n$ is a homeomorphism
and for $1\leq j\leq n$,
$|(i_j,\underline{\mbox{\hspace{.9ex}}})|=|X
(i_j,\underline{\mbox{\hspace{.9ex}}})|= Y(i_j)$. Hence, the map
\[
p_n=\langle Y(i_1),\ldots,Y(i_n)\rangle= \langle
|k^1|,\ldots,|k^n|\rangle \circ  |\langle
(i_1,\underline{\mbox{\hspace{.9ex}}}),\ldots,(i_n,\underline{\mbox{\hspace{.9ex}}})\rangle|,
\]
as a composition of a homeomorphism with $|p_{n
\underline{\mbox{\hspace{.9ex}}}}|$, is a homotopy equivalence
between $Y_n$ and $(Y_1)^n$. We conclude that $Y$ is Segal's, and
by Proposition~3.1, $Y_1$ is a homotopy associative H-space.

If a simplicial space is Segal's, then its realization is
path-connected. This is because its value at $[0]$ is contractible
and therefore path-connected (see \cite[Lemma~11.11]{M72}). Since
$X_{1\,\underline{\mbox{\hspace{.9ex}}}}$ is Segal's, we conclude
that $Y_1$ is path-connected. Moreover, it is grouplike since
every path-connected homotopy associative H-space, which is a
CW-complex, is grouplike (see \cite[Proposition~8.4.4]{A11}).

Applying Proposition~3.2 to $Y$, we obtain that $Y_1\simeq
\Omega|Y|$. Hence,
\[
X_{11}\simeq \Omega Y_1\simeq \Omega (\Omega|Y|)= \Omega^2|X|.
\]

\vspace{-2.2em}

\mbox{\hspace{2em}}\qed


\begin{thebibliography}{99} 

\bibitem{A11} {\sc M.\ Arkowitz}, \textbf{\textit{Introduction to Homotopy Theory}}, Springer, Berlin, 2011


\bibitem{GZ67} {\sc P.\ Gabriel} and {\sc M.\ Zisman}, \textbf{\textit{Calculus of Fractions and Homotopy Theory}}, Ergebnisse der Mathematik und Ihrer Grenzgebiete,
vol.\ 35, Springer, Berlin, 1967


\bibitem{ML71} {\sc S.\ Mac Lane}, \textbf{\textit{Categories for the Working
Mathematician}}, Springer, Berlin, 1971 (expanded second edition,
1998)

\bibitem{M67} {\sc J.P.\ May},
\textbf{\textit{Simplicial Objects in Algebraic Topology}}, The
University of Chicago Press, Chicago, 1967

\bibitem{M72} --------, \textbf{\textit{The Geometry of Iterated Loop Spaces}}, Lecture Notes in Mathematics, vol.\ 271, Springer,
Berlin, 1972

\bibitem{M74} --------, {\it $E_\infty$-spaces, group completions and permutative
categories}, \textbf{\textit{New Developments in Topology}} (G.\
Segal, editor), London Mathematical Society Lecture Notes Series,
vol.\ 11, Cambridge University Press, 1974, pp.\ 153-231


\bibitem{PT13} {\sc Z.\ Petri\' c} and {\sc T.\ Trimble},
{\it Symmetric bimonoidal intermuting categories and
$\omega\times\omega$ reduced bar constructions},
\textbf{\textit{Applied Categorical Structures}}, vol.\ 22 (2014),
pp.\ 467-499 (arXiv:0906.2954)

\bibitem{P14} {\sc Z.\ Petri\' c},
{\it Segal's multisimplicial spaces}, \textbf{\textit{Publications
de l' Institut Mathematique}}, tome 97 (111) (2015), pp.\ 11-21
(arXiv:1407.3914)

\bibitem{S74} {\sc G.\ Segal}, {\it Categories and cohomology theories},
\textbf{\textit{Topology}}, vol.\ 13 (1974), pp.\ 293-312

\bibitem{T79} {\sc R.W.\ Thomason}, {\it Homotopy colimits in the category of small categories},
\textbf{\textit{Mathematical Proceedings of the Cambridge
Philosophical Society}}, vol.\ 85, 91 (1979), pp.\ 91-109

\end{thebibliography}
\end{document}